\documentstyle[12pt,epsf]{article}

\font\twelvesym=msbm10 at 12pt

\font\sevensym=msbm7
\font\fivesym=msbm5

\newfam\symfam
\textfont\symfam=\twelvesym
\scriptfont\symfam=\sevensym
\scriptscriptfont\symfam=\fivesym

\font\twelvegoth=eufm10 at 12pt
\font\tengoth=eufm10
\font\sevengoth=eufm7
\font\fivegoth=eufm5

\newfam\gothfam
\textfont\gothfam=\twelvegoth
\scriptfont\gothfam=\sevengoth
\scriptscriptfont\gothfam=\fivegoth
\def\goth{\fam\gothfam\tengoth}

\def\SG{{\goth S}}

\def\Des{{\rm Des\,}}

\newtheorem{example}{Example}[section]

\newtheorem{theorem}[example]{Theorem}

\newtheorem{corollary}[example]{Corollary}

\def\qed{\hspace{3.5mm} \hfill \vbox{\hrule height 3pt depth 2 pt width 2mm}}

\def\U0{{\cal U}_0(gl_N)}

\def\<{\langle}
\def\>{\rangle}

    \newdimen\squaresize \squaresize=18pt
\newdimen\thickness \thickness=0.5pt
 
\def\square#1{\hbox{\vrule width \thickness
   \vbox to \squaresize{\hrule height \thickness\vss
      \hbox to \squaresize{\hss#1\hss}
   \vss\hrule height\thickness}
\unskip\vrule width \thickness}
\kern-\thickness}
 
\def\vsquare#1{\vbox{\square{$#1$}}\kern-\thickness}
\def\blank{\omit\hskip\squaresize}
 
\def\young#1{
\vbox{\smallskip\offinterlineskip
\halign{&\vsquare{##}\cr #1}}}

\title{\bf 
The Cycle Enumerator of Unimodal Permutations
}

\author{
Jean--Yves {\sc Thibon}
\thanks{\  IGM -- Universit\'e de Marne--la--Vall\'ee -- 
Cit\'e Descartes --- 5 Boulevard Descartes --- Champs-sur-Marne --
77454 Marne-la-Vall\'ee cedex 2 -- France --
{\it e-mail}{\,}:
{\tt jyt@univ-mlv.fr}}
\\
\hbox{}
}

\date{}

\begin{document}

\maketitle

\footnotesize\noindent
{\bf Abstract:  } We give a generating function for the number of unimodal
permutations with a given cycle structure.

\normalsize

\section{Introduction}

The problem of determining the number of unimodal permutations
with a given cycle structure was proposed by Rogers \cite{Ro},
motivated by problems related to the structure of periodic orbits of
one dimensional dynamical systems. Partial results
have since been obtained by Weiss and Rogers \cite{WR} and
more recently by Gannon \cite{Ga}. The proofs of these
results proceed by construction of explicit bijections. In this
note, we propose another approach,
leading to a complete solution in the form of a 
generating function. Our main tool is  the theory of symmetric
functions.

Recall that a permutation $\sigma\in\SG_n$ is said to be {\em unimodal}
if for some $m\in I_n=\{1,2,\ldots,n\}$,
$$
\sigma(1)<\sigma(2)<\ldots<\sigma(m)>\sigma(m+1)>\ldots >\sigma(n)\,.
$$
The set of such permutations will be denoted by $U_n$. For a partition
$\alpha=(1^{a_1}2^{a_2}\cdots n^{a_n})$ of $n$, let $u_\alpha$ be the number
of unimodal permutations  having $\alpha$ as cycle type (i.e., 
having $a_1$ fixed points, $a_2$ 2-cycles, etc.), and let
$p_\alpha=p_1^{a_1}p_2^{a_2}\cdots p_n^{a_n}$ (here, the $p_k$, $k\ge 1$
are independent variables, interpreted in the sequel as power-sum
symmetric functions).

\begin{theorem}\label{th:1}
 Let $c_n=\displaystyle {1\over n}\sum_{d|n,\,d\ {\rm odd}}\mu(d)2^{n/d-1}$,
where $\mu$ is the M\"obius function. The cycle enumerator of
unimodal permutations is given by
$$
1+2\sum_{|\alpha|\ge 1} u_\alpha p_\alpha = \prod_{k\ge 1}
\left( {1+p_k\over 1-p_k}\right)^{c_k}\,.
$$
\end{theorem}

The coefficient of $p_n$ being obviously equal to $2c_n$,
we recover the result of Weiss and Rogers \cite{WR}:

\begin{corollary}
The number of transitive unimodal permutations in $\SG_n$ is equal
to $c_n$.
\end{corollary}

The specialization $p_n=t^n$ leads to various generating functions.
The simplest example is of course the identity
$$
\prod_{k\ge 1}
\left( {1+t^k\over 1-t^k}\right)^{c_k} = {1\over 1-2t}
$$
which is easily checked by expanding the logarithms of both sides,
and amounts to the well-known fact that $|U_n|=2^{n-1}$.
Then, if we denote by $u_n^{(k)}$ the number of $\sigma\in U_n$
having no $k$-cycle, we find

\begin{corollary}
$$
1+2\sum_{n\ge 1} u_n^{(k)}t^n = \left({1-t^k\over 1+t^k}\right)^{c_k}{1\over 1-2t}\,.
$$
\end{corollary}
For example, the number of unimodal permutations without fixed points
(derangements) is seen to be $u_n^{(1)}={1\over 3}(2^{n-1}+(-1)^n)$,
in accordance with \cite{Ga}.

Finally, we note that the generating function for the number $v_n^{(m)}$
of unimodal permutations satifying $\sigma^m=1$ is a rational
function:

\begin{corollary}
$$
1+2\sum_{n\ge 1} v_n^{(m)}t^n = \prod_{d|m}\left({1+t^d\over 1-t^d}\right)^{c_d}\,.
$$
\end{corollary}

\section{Background}

An integer $d\in\{1,2,\ldots,n-1\}$ is said to be a {\em descent} of
a permutation $\sigma\in\SG_n$ if $\sigma(i)>\sigma(i+1)$.
Let $\Des(\sigma)=D=\{d_1,d_2,\ldots,d_r\}$ be the descent set
of $\sigma$. It is customary to encode it by a composition
$I=(i_1,\ldots,i_{r+1})$ of $n$ in $r+1$ parts, called the
{\em shape} of $\sigma$, defined by $i_1=d_1$, $i_2=d_2-d_1$,
$\ldots$, $i_{r+1}=n-d_r$. 
We write $I=C(D)=C(\sigma)$ and $D=\Des(I)$.
These definitions are best visualized
by writing $\sigma$ as a ribbon diagram, e.g., with 
$\sigma=42687135$,
$$
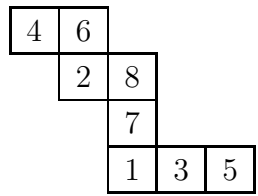
\young{4&6&\blank &\blank &\blank \cr
       \blank &2&8&\blank &\blank \cr
       \blank &\blank & 7 &\blank & \blank \cr
       \blank &\blank &1&3&5\cr}
$$
the lengths of the rows form the composition $(2,2,1,3)$ of $8$,
corresponding to the descent set $\{2,4,5\}$.

Thus, unimodal permutations are those of shape $(n-k,1^k)$
for some $k=0,\ldots,n-1$. The number
$a_{\alpha,I}$ of permutations of shape $I$ and cycle type $\alpha$
is known to be equal to a scalar product of symmetric functions,
which may sometimes be evaluated in closed form \cite{GR}.

Let $Sym=Sym(X)$ be the ring of symmetric functions in the
set of indeterminates $X=\{x_1,x_2,\ldots\}$, endowed with its
standard scalar product for which the Schur functions $s_\lambda$
form an orthonormal basis (see, e.g., \cite{Mcd}). 

Let $h_n=h_n(X)=s_{(n)}(X)$ be the complete homogeneous symmetric functions,
defined by the generating series
$$
\sigma_t(X)=\prod_{i\ge 1}(1-tx_i)^{-1}=\sum_{n\ge 0}h_nt^n
$$
and $p_n=\displaystyle\sum_{i\ge 1}x_i^n$ be the power sums.
We shall need the Witt symmetric functions
$$
\ell_n={1\over n}\sum_{d|n}\mu(d)p_d^{n/d}
$$
and
$$
L_\alpha=(h_{a_1}\circ \ell_1)(h_{a_2}\circ \ell_2)\cdots (h_{a_m}\circ\ell_n)
$$
where $g\circ f$ denotes the {\em plethysm} of $f$ by $g$,
and the MacMahon symmetric functions (or ribbon Schur functions,
denoted by $h_I$ in \cite{MM})
$$
r_I=
\left|\matrix{
h_{i_1} & h_{i_1+i_2} & \cdots & h_{i_1+i_2+\cdots+i_m} \cr
1       & h_{i_2}     & \cdots & h_{i_2+\cdots+i_m} \cr 
0       & 1           &\cdots  &  h_{i_3+\cdots+i_m} \cr
\vdots  & \vdots      & \ddots & \vdots\cr
0       &  0          & \cdots & h_{i_m}\cr
}\right|
$$
defined for any composition $I=(i_1,i_2,\ldots,i_m)$.

\begin{theorem}[{\rm Gessel-Reutenauer \cite{GR}}]\label{th:GR} 
The number $a_{\alpha,I}$ of permutations
of shape $I$ and cycle type $\alpha$ is equal to $\<r_I,L_\alpha\>$.
\end{theorem}

This result was originally obtained by means of a  combinatorial argument.
Recently, J\"ollenbeck and Reutenauer have found a simple algebraic proof,
relying on the following formula, which we will also need in the sequel:

\begin{theorem}[{\rm Scharf-Thibon \cite{ST}}]\label{th:ST}
Let $X$ and $Y$ be two independent sets
of indeterminates, and set ${\cal L}(X,Y)=\sum_\alpha p_\alpha(X)L_\alpha(Y)$.
Then,
$$ {\cal L}(X,Y)={\cal L}(Y,X)\,.
$$
\end{theorem}

Theorem \ref{th:1} is a simple consequence of the above results.
The easiest way to derive it is to make use of a
classical ``$\lambda$-ring'' trick for dealing with hook Schur
functions. This works as follows. The argument $X$ of a symmetric
function is identified with the {\em formal} sum of its elements:
$X=x_1+x_2+\cdots$. Then, $X\cup Y$ is identified with $X+Y$,
and one can define $XY=\sum_{i,j}x_iy_j$. The power-sums behave
then as ring homomorphisms: $p_n(X+Y)=p_n(X)+p_n(Y)$ and $p_n(XY)=p_n(X)p_n(Y)$.
Let $q$ be another indeterminate. Under the specialization $x_k=q^{k-1}$,
$X$ becomes identified with the power series $(1-q)^{-1}$, with
formal inverse $1-q$. The symmetric functions of $1-q$ are then specified
by the condition $p_n(1-q)=1-q^n$. 

With this notation, we can state the well-known identity
\begin{equation}\label{eq:1}
h_n((1-q)X)=(1-q)\sum_{k=0}^{n-1}(-q)^ks_{(n-k,1^k)}(X)\,.
\end{equation}

\section{Proof of Theorem \ref{th:1}}

Since $r_{(n-k,1^k)}=s_{(n-k,1^k)}$, we have by Theorem \ref{th:GR}
and (\ref{eq:1})
$$
u_\alpha=\sum_{k=0}^{n-1}\<s_{(n-k,1^k)},L_\alpha\>
=\left\< {h_n((1-q)X)\over 1-q}\,,\, L_\alpha(X)\right\>_{q=-1}
$$
$$
={1\over 2}\<h_n(X)\,,\,L_\alpha((1-q)X)\>_{q=-1}={1\over 2}L_\alpha(1-q)|_{q=-1}\,.$$
Therefore,
$$
1+2\sum_{|\alpha|\ge 1} u_\alpha p_\alpha(X) =
\sum_{|\alpha|\ge 0} L_\alpha(1-q)|_{q=-1}p_\alpha(X)
=\sum_{|\alpha|\ge 0}p_\alpha(1-q)|_{q=-1}L_\alpha(X)
$$
(by Theorem \ref{th:ST}).
But $p_\alpha(1-q)|_{q=-1}$ is zero as soon as $\alpha$ has an even
part, and is equal to $2^{l(\alpha)}$ otherwise, where $l(\alpha)$ is
the length of $\alpha$. Hence,
$$
1+2\sum_{|\alpha|\ge 1} u_\alpha p_\alpha(X) 
=\sum_{\alpha \ {\rm odd}} 2^{l(\alpha)}L_\alpha(X)=\sigma_2\circ
\left[\sum_{k\ge 0} \ell_{2k+1}(X)\right]
$$
$$
=\exp\left\{
\sum_{j\ge 1} {2^jp_j\over j}\circ\left[
\sum_{k\ge 0} {1\over 2k+1} \sum_{d|k}\mu(d)p_d^{(2k+1)/d}\right]\right\}
$$
$$
=\exp\left\{
\sum_{j\ge 1}{2^jp_j\over j}\circ\left[ 
\sum_{a\ge 0} {\mu(2a+1)\over 2a+1}\sum_{b\ge 0}{p_{2a+1}^b\over 2b+1}
\right]\right\}
$$
$$
=\exp\left\{ 
\sum_{j\ge 1}\sum_{a\ge 0}
{2^j\mu(2a+1)\over (2a+1)j}\log\left({1+p_{(2a+1)j}\over 1-p_{(2a+1)j}}\right)
\right\}
$$
$$
=\prod_{j\ge 1}\prod_{a\ge 0}
\left(
{1+p_{(2a+1)j}\over 1-p_{(2a+1)j}}\right)^{{\mu(2a+1)2^{j-1}\over (2a+1)j}}
=\prod_{n\ge 1} \left( {1+p_n\over 1-p_n}\right)^{c_n}
$$
\qed

\section{Further comments}

\subsection{A $q$-analogue of Theorem \ref{th:1}}

The previous calculation can be carried out without setting
$q=-1$, and this yields the enumeration of unimodal permutations
according to  cycle type and position of the maximum.
Setting $u_\alpha(q)=\sum_k u_{\alpha,k}q^k$ where $u_{\alpha,k}=\<s_{(n-k,1^k)},L_\alpha\>$
is the number of elements of $U_n$ with cycle type $\alpha$ and maximum
at $m=n-k$, we find
$$
1+(1-q)\sum_{n\ge 1}\sum_{|\alpha|=n}\sum_{k=0}^{n-1}(-q)^ku_{\alpha,k}p_\alpha(X)
=
\<\sigma_1((1-q)Y)\,,\,{\cal L}(Y,X)\>_Y
$$
$$
={\cal L}(1-q,X) =\prod_{k\ge 1}\sigma_{1-q^k}\circ\ell_k(X)
$$
$$
=\exp\left\{
\sum_{k\ge 1}\sum_{m\ge 1}{(1-q^k)^m\over m} p_m\circ
\left[
{1\over k}\sum_{d|k}\mu(d)p_d^{k/d}
\right]
\right\}
$$
$$
=\exp\left\{ 
\sum_{n\ge 1}\sum_{a\ge 1} \ell_n(1-q^a) {p_n^a\over a}
\right\}
$$
(after some simple transformations), and setting
\begin{equation}
\ell_n(1-q)=
{1\over n}\sum_{d|n}\mu(d)(1-q^d)^{n/d}
=
\sum_k a_{nk}q^k
\end{equation}
we obtain
\begin{theorem}\label{th:q}
$$
{\cal L}(1-q,X)=\prod_{n\ge 1}\prod_k\left(1-q^kp_n\right)^{-a_{nk}}\,.
$$
\end{theorem}
The first values of the $a_{nk}$ are given by
$$
\ell_1(1-q)=1-q,\ 
\ell_2(1-q)=-q+q^2,\
\ell_3(1-q)=-q+q^2,\
\ell_4(1-q)=-q+2q^2-q^3,\ldots
$$ 
so that
$$
{\cal L}(1-q,X)=
\left({1-qp_1 \over 1- p_1}\right)
\left({1-qp_2 \over 1- q^2p_2}\right)
\left({1-qp_1 \over 1- q^2p_3}\right) 
{(1-qp_4)(1-q^3p_4)\over (1-q^2p_4)^2}\cdots
$$
and
$$
\sum_{|\alpha|\ge 1}u_\alpha(q)=
{1\over 1-q}\left[{\cal L}(1-q,X)-1\right]_{q\rightarrow -q}
$$
$$
=p_1 + p_{11}+qp_2 + p_{111}+(q+q^2)p_{21}+qp_3
$$
$$
+p_{1111}+(q+q^2)p_{211}+q^2p_{22}+(q+q^2)p_{31}+(q+q^2)p_4+\cdots
$$

\subsection{Unimodal permutations and iterated brackets}

One way to generate the set of unimodal permutations is to
expand the iterated Lie bracket
\begin{equation}\label{eq:bra}
[x_1,[x_2,[x_3,\cdots , [x_{n-1},x_n] \cdots ]]]
= \sum_{\sigma\in U_n} 
\nu(\sigma)
x_{\sigma(1)}x_{\sigma(2)}\cdots x_{\sigma(n)}\,,
\end{equation}
where $\nu(\sigma)=(-1)^k$ if the descent composition of $\sigma$
is $(n-k,1^k)$.
Clearly, replacing the bracket by its $q$-analogue $[a,b]_q=ab-qba$ leads to
the same formula, with $\nu$ replaced by $\nu_q(\sigma)=(-q)^k$.

In \cite{BL}, Blessenohl and Laue prove the following identity:
\begin{equation}\label{eq:BL}
\sum_{\sigma\in U_n\cap C_\alpha} \nu(\sigma)
=
\left\{\matrix{
\alpha(d) & \ {\rm if} \ \alpha=(d^m),\cr
0      & {\rm otherwise,}\cr}
\right.
\end{equation}
where $C_\alpha$ denotes the set of permutations of cycle type $\alpha$.
This result can be easily derived from a calculation of symmetric
functions. Indeed,
$$
\sum_{\sigma\in U_n\cap C_\alpha} \nu(\sigma) 
=
\sum_{k=0}^{n-1}\<s_{(n-k,1^k)},L_\alpha\>=\<p_n,L_\alpha\>\,.
$$
This scalar product is non zero only when the power-sum expansion
of $L_\alpha$ contains a term proportional to $p_n$, which is the case
only when the product $L_\alpha=\prod_i h_{a_i}[\ell_i]$ is reduced
to a single term 
$$
h_a\circ\ell_d=\sum_{|\alpha|=a}{1\over z_\alpha}p_\alpha\circ
\left[
{1\over d}\sum_{j|d}\mu(j)p_j^{d/j}\right]
$$
where the term involving $p_n$ is
$$
{p_a\over a}\circ\left[ {\mu(d)\over d}p_d\right]
=\mu(d){p_n\over n}\,,
$$
whence the result.

The most natural $q$-analogue of (\ref{eq:BL}) would be
the evaluation of
$$
\sum_{\sigma\in U_n\cap C_\alpha} \nu_q(\sigma)={1\over 1-q}L_\alpha(1-q)\,.
$$
By Theorem \ref{th:ST}, we see that this is precisely the content
of Theorem \ref{th:q}, for
 $${\cal L}(1-q,X)={\cal L}(X,1-q)=\sum_{\alpha}L_\alpha(1-q)p_\alpha(X)\,.
$$

\subsection{Products of unimodal permutations}

The properties of unimodal permutations are directly related
to quasi-symmetric functions \cite{Ge} and noncommutative
symmetric functions \cite{NCSF1,NCSF2}. To give one more
example, let us remark that the result of Kreweras \cite{Kr}
that any unimodal permutation of $\SG_n$ is, in exactly $n$
ways, the product of two unimodal permutations, appears now
as a special case of Corollary 5.9 of \cite{NCSF2}, describing
the product $U_n(x)U_n(y)$ in the group algebra, where
$$
U_n(q)=\sum_{k=0}^{n-1}q^k\sum_{C(\sigma)=(n-k,1^k)}\sigma\,.
$$
The coefficient of $\sigma$ in the product is
$y^r[r+1]_{x/y} + x^{r+1}[n-r-1]_{xy}$ if $C(\sigma)=(n-r,1^r)$
(unimodal case), $x^{r+s-1}y^{t+s-1}(y+1)$ if
$C(\sigma)=(n-r-s-t,1^t,s,1^r)$ (bimodal case), and $0$ otherwise
(here, $[n]_q=(1-q^n)/(1-q)$).

\subsection{Grassmannian permutations}

All the calculations in this paper could be easily adapted
to treat the case of Grassmannian permutations, i.e. permutations
with at most one descent. The relevant specialization
would be here $X=1+q$.


\end{document}